\documentclass[11pt, reqno]{amsart}

\usepackage{amssymb}
\usepackage{amsmath}
\usepackage{mathrsfs}
\usepackage{amsfonts}
\usepackage{color}
\usepackage{vmargin}
\usepackage{amsthm}
\usepackage{graphicx}

\usepackage{amssymb}
\usepackage{amsmath}
\usepackage{enumerate}
\usepackage{mathrsfs}
\usepackage{amsfonts}
\usepackage{color}
\usepackage{vmargin}
\usepackage{amsthm}
\usepackage{graphicx} 
\usepackage{esint}

\usepackage{latexsym}
\usepackage{xcolor}
\usepackage{mleftright}
\usepackage{tikz}
\usepackage{tikz-3dplot}
\usetikzlibrary{calc}
\usetikzlibrary{arrows}
\usetikzlibrary{positioning}
\usepackage{mathabx}
\usepackage{bbm}
\usepackage[hidelinks]{hyperref}
\usepackage{mathrsfs}
\usepackage{enumitem}

\setmarginsrb{20mm}{20mm}{20mm}{20mm}{10mm}{10mm}{10mm}{10mm}

\RequirePackage{amsthm}

\theoremstyle{definition}
\newtheorem{theorem}{Theorem}
\numberwithin{theorem}{section}

\theoremstyle{definition}
\newtheorem{remark}[theorem]{Remark}
 
 \theoremstyle{definition}
\newtheorem{lemma}[theorem]{Lemma}

\theoremstyle{definition}
\newtheorem{sublemma}{Sublemma}[theorem]

\theoremstyle{definition}
\newtheorem{definition}[theorem]{Definition}

\RequirePackage{amsthm}
\newcommand{\brac}[1]{\left({#1}\right)}

\theoremstyle{definition}
\newtheorem{claim}{Claim}
\numberwithin{claim}{section}

\theoremstyle{definition}
\newtheorem{proposition}[theorem]{Proposition}

\theoremstyle{definition}
\newtheorem{corollary}[theorem]{Corollary}

\theoremstyle{definition}

\theoremstyle{definition}

\theoremstyle{definition}
\newtheorem*{theorem*}{Theorem}

\newcommand{\PS}{\mathbf{PS}(\delta)}
\newcommand{\MS}{\mathbf{MS}(\delta)}
\newcommand{\PN}{\mathbf{PN}(\delta)}
\newcommand{\MN}{\mathbf{MN}(\delta)}

\newcommand{\PSE}{\mathbf{PSE}(\delta)}
\newcommand{\MSE}{\mathbf{MSE}(\delta)}
\newcommand{\PNE}{\mathbf{PNE}(\delta)}
\newcommand{\MNE}{\mathbf{MNE}(\delta)}

\setmarginsrb{20mm}{20mm}{20mm}{20mm}{10mm}{10mm}{10mm}{10mm}

\newtheoremstyle{remark}
  {}{}{}{}{\bfseries}{.}{.5em}{{\thmname{#1 }}{\thmnumber{#2}}{\thmnote{ (#3)}}}
\theoremstyle{remboldstyle}

\begin{document}

\title[Maximal $\delta$-separated sets and \textbf{AC}]{Maximal $\delta$-separated sets in separable metric spaces and weak forms of choice}

\author[M. Dybowski]{Micha\l{} Dybowski}
\address{Faculty of Mathematics and Information Scineces,
Warsaw University of Technology,
Pl. Politechniki 1, 00-661 Warsaw, Poland}
\email{michal.dybowski.dokt@pw.edu.pl}

\author[P. G\'{o}rka]{Przemyslaw G\'{o}rka}
\address{Faculty of Mathematics and Information Scineces,
Warsaw University of Technology,
Pl. Politechniki 1, 00-661 Warsaw, Poland}
\email{przemyslaw.gorka@pw.edu.pl}

\author[P. Howard]{Paul Howard}
\address{2770 Ember Way, Ann Arbor, MI 48104, USA}
\email{phoward@emich.edu}


\begin{abstract}
We show that the statement ``In every separable pseudometric space there is a maximal non-strictly $\delta$-separated set.'' implies the axiom of choice for countable families of sets.  This gives answers to a question of Dybowski and G\'{o}rka \cite{Dyb-Gorka1}.  We also prove several related results.
\end{abstract} 

\maketitle

\section{Introduction}
We will use the following standard notation:  \textbf{ZF} for Zermelo-Fraenkel set theory, \textbf{ZFC} for Zermelo-Fraenkel set theory with the axiom of choice, \textbf{ZFA} for \textbf{ZF} weakened to permit the existence of atoms and \textbf{ZFA + AC} for \textbf{ZFA} with the addition of the axiom of choice.

All of our positive results about the implications between various consequences of the axiom of choice are theorems of \textbf{ZFA} and therefore theorems of the stronger theory \textbf{ZF}.  Our independence results all hold in the theory \textbf{ZF} although many of them are  first proved in \textbf{ZFA} and then transferred to \textbf{ZF} using a transfer theorem of Pincus.
 
Let $\delta>0$, the subject of this paper is the deductive strength of the following theorem of \textbf{ZFA + AC}:
\begin{equation*}
\begin{aligned}
&\mbox{\it For every separable pseudometric space $(X, d)$ },\\ & X \mbox{\it has a maximal non-strictly } \delta\mbox{\it -separated set.}
\end{aligned}
\end{equation*}
 and three closely related theorems (see Definition \ref{D:Weakchoice}, (\ref{D:WeakchoiceSps}) - (\ref{D:WeakchoiceSmn})).
Our abbreviation for this theorem of $\mathbf{ZFA + AC}$ is $\PN$ (\textbf{P} for `pseudometric' and \textbf{N} for `non-stricty'). 

We begin with the definitions we will require.
\begin{definition}
(Metric space and set theoretic terminology)
\begin{enumerate}
\item We use $\omega$ for the set of natural numbers.  That is, $\omega = \{0, 1, 2, \ldots \}$.
\item $\mathbb{Z}$ denotes the set of integers and $\mathbb{Z}^+ = \{ 1, 2, 3, \ldots \}$ is the set of positive integers.
\item 
A \emph{metric space} is a pair $(X,d)$ where $X$ is a set and $d: X \times X \to [0, \infty)$ is a function which satisfies: For all $x$, $y$ and $z$ in $X$
 \begin{enumerate}
 \item \label{I:dxx=0} $d(x,y) = 0 \iff x=y$ 
 \item \label{I:dxy=dyx} $d(x,y) = d(y,x)$
 \item \label{I:triangle} $d(x,z) \le d(x,y) + d(y,z)$
 \end{enumerate}
\item 
 A \emph{pseudometric space} is a pair $(X,d)$ satisfying items (\ref{I:dxy=dyx}), (\ref{I:triangle}) above and for every $x \in X$ $d(x,x)=0$.
 \item 
An \emph{ultrametric space} is a metric space $(X,d)$  such that  for every $x, y, z \in X$
\[
d(x,z) \le \mbox{max}(d(x,y), d(y,z)).
\]
\item Assume that $(X,d)$ is either a metric space or a pseudometric space, that $Y  \subseteq X$ and that $\delta > 0$ is a real number. Then
 \begin{enumerate}
 \item $Y$ is a \emph{strictly $\delta$-separated set} (in $(X,d)$) if for all distinct points $x$ and $y$ in $Y$, $d(x,y) > \delta$.
 \item $Y$ is a \emph{non-strictly $\delta$-separated set} (in $(X,d)$) if for all distinct points $x$ and $y$ in $Y$, $d(x,y) \ge \delta$.
 \end{enumerate}
\end{enumerate} 
\end{definition}

\begin{definition} \label{D:Weakchoice}
(Weak choice axioms)
\begin{enumerate}
\item \textbf{AC} (Form 1 in \cite{HR}) is the axiom of choice:  For every set $X$ of non-empty sets, there is a function $f : X \to \bigcup X$ such that $\forall y \in X, f(y) \in y$.
\item {\textbf{CC}} (Form 8 in \cite{HR}): For every countable set $X$ of non-empty sets, there is a function $f : X \to \bigcup X$ such that $\forall y \in X, f(y) \in y$. 
\item $\mathbf{CC}(\mathbb{R})$ (Form 94 in \cite{HR}): For every countable set $X$ of non-empty sets of real numbers, there is a function $f : X \to \bigcup X$ such that $\forall y \in X, f(y) \in y$.
\item A relation $R$ on set $X$ is \emph{entire (on $X$)} if for every $x \in X$ there exists $y \in X$ such that $x \mathrel{R} y$.
\item {\textbf{DC}} (Form 43 in \cite{HR}) is ``For every non-empty set $X$ and every entire relation $R$ on $X$ there is a sequence $\{ x_i \}_{i = 1} ^\infty \subset X$ such that $x_i \mathrel{R} x_{i+1}$ for all $i \in \mathbb{Z}^+$.
\vskip.1in
\item[] Each of the following eight statements has a parameter $\delta$, where $\delta$ is a positive real number.
\vskip.1in
\item \label{D:WeakchoiceSps} $\PS $:  For every separable pseudometric space $(X,d)$, $X$ has a maximal strictly $\delta$-separated subset.
\item $\MS$:  For every separable metric space $(X,d)$, $X$ has a maximal strictly $\delta$-separated subset.
\item $\PN$:  For every separable pseudometric space $(X,d)$, $X$ has a maximal non-strictly $\delta$-separated subset.
\item \label{D:WeakchoiceSmn} $\MN$:  For every separable metric space $(X,d)$, $X$ has a maximal non-strictly $\delta$-separated subset.

\item \label{D:WeakchoicePS} $\PSE $:  For every separable pseudometric space $(X,d)$, if $S \subseteq X$ is strictly $\delta$-separated set then there is an $S' \subseteq X$ such that $S \subseteq S'$ and $S'$ is a maximal strictly $\delta$-separated set.
\item $\MSE$:  For every separable metric space $(X,d)$, if $S \subseteq X$ is strictly $\delta$-separated set then there is an $S' \subseteq X$ such that $S \subseteq S'$ and $S'$ is a maximal strictly $\delta$-separated set.
\item $\PNE$:  For every separable pseudometric space $(X,d)$, if $S \subseteq X$ is non-strictly $\delta$-separated set then there is an $S' \subseteq X$ such that $S \subseteq S'$ and $S'$ is a maximal non-strictly $\delta$-separated set.
\item \label{D:WeakchoiceMN} $\MNE$:  For every separable metric space $(X,d)$, if $S \subseteq X$ is non-strictly $\delta$-separated set then there is an $S' \subseteq X$ such that $S \subseteq S'$ and $S'$ is a maximal non-strictly $\delta$-separated set.

\end{enumerate}
\end{definition}
\par 
In \cite{Dyb-Gorka1} Dybowski and G\'{o}rka studied $\PS$ and $\MS$.  They proved that several strengthenings of $\PS$ are equivalent to \textbf{AC} in \textbf{ZF}.  To be more specific, in their Proposition 4.2 they proved that the following statement is equivalent to \textbf{AC}:
\begin{equation}
\begin{aligned}
\label{E:mspaceContainsDelSepSet}
&\mbox{\it For every pseudometric space $(X, d)$, and every $\delta >0$} \\ & X \mbox{\it has a maximal strictly } \delta\mbox{\it -separated set.}
\end{aligned}
\end{equation}
(Note that the ``separable'' requirement has been removed.)
\par
They also showed that if ``pseudometric'' is replaced by ``metric'' or by ``ultrametric'' in (\ref{E:mspaceContainsDelSepSet}) each of the resulting statements is equivalent to \textbf{AC}. 
\begin{remark}[\textbf{ZFA}]\label{rem1}
For all $\delta \in \mathbb{R}^+$, $\mathbf{DC} \implies \PS + \PN + \MN + \MS$.
\end{remark}
\begin{proof}
The implication $\mathbf{DC} \implies \PS$ was proved before (see \cite[Corollary 4.5]{Dyb-Gorka1}). The proof given can be modified to show that \textbf{DC} implies $\PN$.
Since $\PN$ implies $\MN$ and $\PS$ implies $\MS$, it follows that \textbf{DC} also implies $\MN$ and $\MS$. 
\end{proof}
Finally they ask in Problem 4.6 about what else can be said concerning the strength, in set theory without \textbf{AC}, of the statements $\PS$ and $\MS$.
\par 
The purpose of this paper is to provide some answers to this question and to related questions involving items (\ref{D:WeakchoiceSps}) through (\ref{D:WeakchoiceMN}) in Definition \ref{D:Weakchoice} above.  Our results are summarized in the following diagram.  (See the paragraph following the proof of Theorem \ref{T:DeltaDoesntMatter} for a description of the abbreviations used in the diagram.)

We use a variation of the notation of \cite{bentley}, namely, if $\mathbf{\Psi}$ and $\mathbf{\Phi}$ are sentences and occur in a box together then they are equivalent in \textbf{ZF}.  If $\mathbf{\Psi}$ implies $\mathbf{\Phi}$ in \textbf{ZF} and we don't know whether of not $\mathbf{\Phi}$ implies $\mathbf{\Psi}$ in \textbf{ZF} we use $\mathbf{\Psi} \to \mathbf{\Phi}$. Sentence $\mathbf{\Psi} \mapsto \mathbf{\Phi}$ means that in \textbf{ZF}, $\mathbf{\Psi}$ implies $\mathbf{\Phi}$ but $\mathbf{\Phi}$ does not imply $\mathbf{\Psi}$.  $\mathbf{\Psi} \ \raisebox{-.2ex}{\rotatebox{90}{\scalebox{1}[1.2]{$\bot$}}}\ \mathbf{\Phi}$ means $\mathbf{\Psi}$ does not imply $\mathbf{\Phi}$ in \textbf{ZF} and we don't know whether or not $\mathbf{\Phi}$ implies $\mathbf{\Psi}$.
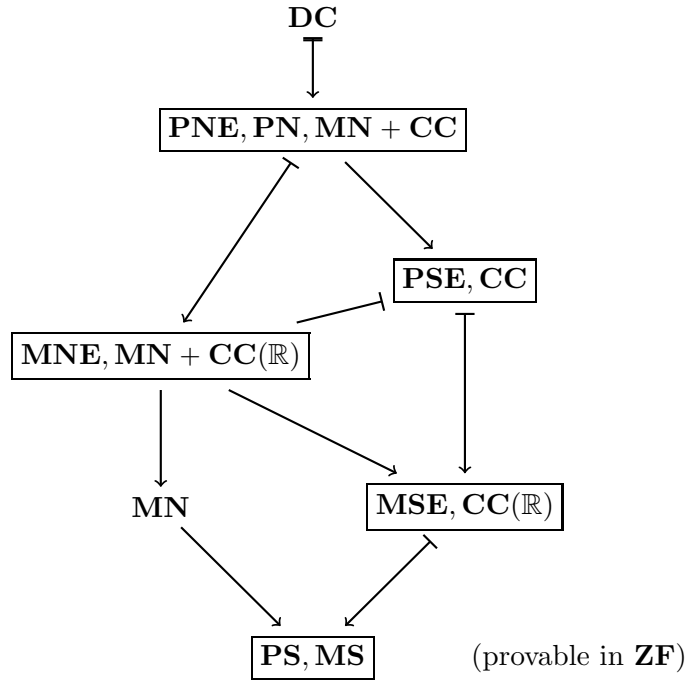
\begin{center}
\begin{tikzpicture}
\tikzset{thick arc/.style={->, black, fill=none, thick, >=stealth, text=black}}
\node (dc) at (0,.5) {$\mathbf{DC}$};
\node (pne_pn) at (0,-1) {$\boxed{\mathbf{PNE,PN,MN+CC}}$};
\node (pse_cc) at (2,-3) {$\boxed{\mathbf{PSE, CC}}$};
\node (mne_mn+ccr) at (-2,-4) {$\boxed{\mathbf{MNE},\mathbf{MN+CC(\mathbb{R})}}$};
\node (mn) at (-2,-6) {$\mathbf{MN}$};
\node (mse_ccr) at (2,-6) {$\boxed{\mathbf{MSE, CC(\mathbb{R})}}$};
\node (ps_ms) at (0,-8) {$\boxed{\mathbf{PS, MS}}$};
\node[right=of ps_ms] (note) {(provable in $\mathbf{ZF}$)};
\draw [thick,{|}{|}->] (dc) -- (pne_pn);
\draw [thick,->] (pne_pn) -- (pse_cc);
\draw [thick,|->] (pne_pn) -- (mne_mn+ccr);
\draw [thick, -|] (mne_mn+ccr) -- (pse_cc);
\draw [thick, ->] (mne_mn+ccr) -- (mn);
\draw [thick, ->] (mne_mn+ccr) -- (mse_ccr);
\draw [thick, |->] (pse_cc) -- (mse_ccr);
\draw [thick, ->] (mn) -- (ps_ms);
\draw [thick, |->] (mse_ccr) -- (ps_ms);
\end{tikzpicture}
\end{center}
\par 
The diagram shows that there are several unsolved problems.  For example
\begin{enumerate}

\item Does $\mathbf{PSE}$ imply $\mathbf{PNE}$?
\item Does \textbf{CC} imply $\mathbf{MN}$ or $\mathbf{MNE}$?
\item Does either of $\mathbf{MN}$ or $\mathbf{MSE}$ imply $\mathbf{MNE}$?
\item \label{Q:ZFtoMN} Is $\mathbf{MN}$ provable in \textbf{ZF}?
\end{enumerate}
We note that if the answer to question (\ref{Q:ZFtoMN}) is in the affirmative\footnote{It is worth to have in mind Theorem \ref{nonstrictZF}.} then all of the other questions are answered and the diagram collapses to
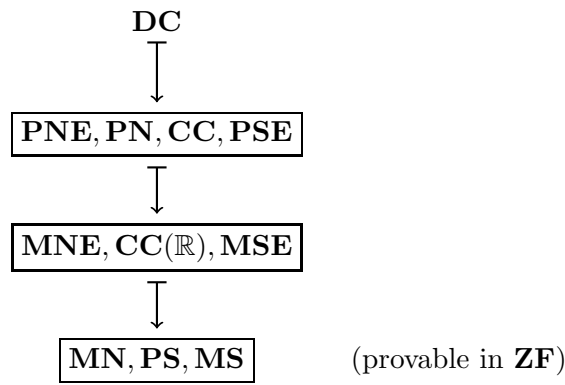
\begin{center}
\end{center}
\begin{center}
\begin{tikzpicture}
\tikzset{thick arc/.style={->, black, fill=none, thick, >=stealth, text=black}}
\node (dc) at (0,0) {$\mathbf{DC}$};
\node (pne_pn_cc_pse) at (0,-1.5) {$\boxed{\mathbf{PNE,PN,CC,PSE}}$};
\node (mne_ccr_mse) at (0,-3) {$\boxed{\mathbf{MNE},\mathbf{CC(\mathbb{R}),MSE}}$};
\node (mn_ps_ms) at (0,-4.5) {$\boxed{\mathbf{MN, PS, MS}}$};
\node[right=of mn_ps_ms] (note) {(provable in $\mathbf{ZF}$)};
\draw [thick,|->] (dc) -- (pne_pn_cc_pse);
\draw [thick,|->] (pne_pn_cc_pse) -- (mne_ccr_mse);
\draw [thick,|->] (mne_ccr_mse) -- (mn_ps_ms);
\end{tikzpicture}
\end{center}

\section{Preliminary Results} \label{S:Prelims}

\begin{proposition}[\textbf{ZFA}] \label{P:XleR} 
If $(X,d)$ is a separable metric space then $|X| \le |\mathbb{R}|$.
\end{proposition}
\begin{proof}
For each point $x \in X$ we can define a (unique) sequence $\hat{x}$ of elements of $D$ (the countable dense set) by letting $\hat{x}_n$ be the first element of $D$ which is in $B(x,1/(n+1))$.  It is clear that the function $x \mapsto \hat{x}$ is one to one so $|X| \le |D^{\aleph_0}| \le |\aleph_0^{\aleph_0}| = |\mathbb{R}|$.
\end{proof}
\begin{theorem}[\textbf{ZFA}]
\label{T:DeltaDoesntMatter}
Assume $\mathbf{\Psi}$ is one of $\mathbf{PS, MS, PN, MN, PSE, MSE, PNE}$ or $\mathbf{MNE}$ and $\delta_1$ and $\delta_2$ are positive real numbers.  Then $\mathbf{\Psi}(\delta_1) \implies \mathbf{\Psi}(\delta_2)$.
\end{theorem}
\begin{proof}
We give an outline of the proof for the case that $\mathbf{\Psi}$ is $\mathbf{PSE}$ and leave the remaining similar proofs for the interested reader.  

Assume $\mathbf{PSE}(\delta_1)$ and let $S$ be a strictly $\delta_2$-separated set in the pseudometric space $(X,p)$. Define a metric $p'$ on $X$ by $p'(x,y) = \frac{ \delta_1}{ \delta_2}p(x,y) $.  Then $(X,p')$ is a pseudometric space and $S$ is a strictly $\delta_1$-separated set in $(X,p')$.  Using $\mathbf{PSE}(\delta_1)$ there is a maximal strictly $\delta_1$-separated set $S'$ in $(X,p')$ such that $S \subset S'$ and therefore $S'$ is a maximal strictly $\delta_2$-separated set in $(X,p)$.
\end{proof}
Having in mind Theorem \ref{T:DeltaDoesntMatter}, for the remainder of the paper we use the following notation. If $\mathbf{\Psi}$ is one of $\mathbf{PS, MS, PN, MN, PSE, MSE, PNE}$ or $\mathbf{MNE}$, we will frequently write $\mathbf{\Psi}$ for any one of the equivalent statements ``$\forall \delta > 0, \mathbf{\Psi}(\delta)$'', ``$\exists \delta > 0, \mathbf{\Psi}(\delta)$'' or, if $\delta$ is a specific positive real number, ``$\mathbf{\Psi}(\delta)$''.

\begin{proposition}[\textbf{ZFA}]\label{CCplusPsi}
Let $(\mathbf{\Psi}_1, \mathbf{\Phi}_1)$ be (\textbf{PS}, \textbf{PSE}) or (\textbf{PN}, \textbf{PNE}) and let $(\mathbf{\Psi}_2, \mathbf{\Phi}_2)$ be (\textbf{MS}, \textbf{MSE}) or (\textbf{MN}, \textbf{MNE}). Then, \textbf{CC} $+$ $\mathbf{\Psi}_1 \implies \mathbf{\Phi}_1 \implies \mathbf{\Psi}_1$ and \textbf{CC}($\mathbb{R}$) $+$ $\mathbf{\Psi}_2 \implies \mathbf{\Phi}_2 \implies \mathbf{\Psi}_2$.
\end{proposition}

\begin{proof}
Let us notice that the implications $\mathbf{\Phi}_1 \implies \mathbf{\Psi}_1$ and $\mathbf{\Phi}_2 \implies \mathbf{\Psi}_2$ are clear. 
We shall prove only the implication \textbf{CC} $+$ $\mathbf{\Psi}_1 \implies \mathbf{\Phi}_1$ 
since \textbf{CC}($\mathbb{R}$) $+$ $\mathbf{\Psi}_2 \implies \mathbf{\Phi}_2$ follows in the similar manner using the fact that all separable metric spaces can be bijectively embedded in $\mathbb{R}$ (see Proposition \ref{P:XleR}). We prove the implication \textbf{CC} $+$ $\mathbf{\Psi}_1 \implies \mathbf{\Phi}_1$  for pair $(\mathbf{\Psi}_1, \mathbf{\Phi}_1) = (\textbf{PN}, \textbf{PNE})$ since the remaining case follows in the similar way.

    Let us assume \textbf{CC} and \textbf{PN} hold. Let $(X, d)$ be a pseudometric separable space, $\delta > 0$ and $S \subset X$ be an arbitrary non-strictly $\delta$-separated set. We define $N_\delta(S) = \bigcup_{x 
    \in S} B\brac{x, \delta}$. Let us consider pseudometric space $\brac{X \setminus N_\delta(S), d}$. Under \textbf{CC} separability is hereditary\footnote{For the implication \textbf{CC}($\mathbb{R}$) $+$ $\mathbf{\Psi}_2 \implies \mathbf{\Phi}_2$ we use analogous fact that for spaces contained in $\mathbb{R}$ separability is hereditary under \textbf{CC}($\mathbb{R}$).} (see \cite[Theorem $1.12$]{bentley}) so $X \setminus N_\delta(S)$ is separable. Thus by virtue of \textbf{PN} space $X \setminus N_\delta(S)$ contains a maximal non-strictly $\delta$-separable set $S_0$. We claim that $S' := S \cup S_0$ is the maximal non-strictly $\delta$-separated set in $X$ containing $S$.

    Indeed, let $x, y \in S'$, $x \neq y$. If $x, y \in S$ or $x, y \in S_0$, then $d(x, y) \ge \delta$. If $x \in S$ and $y \in S_0$, then $y \notin N_\delta(S)$ so $y \notin B(x, \delta)$. It proves that $S'$ is non-strictly $\delta$-separated set.

    Let us suppose that $S'$ is not a maximal non-strictly $\delta$-separated set i.e. there exists $x_0 \notin S'$ such that $d(x_0, x) \ge \delta$ for all $x \in S'$. In particular $d(x_0, x) \ge \delta$ for all $x \in S$ so $x_0 \notin N_\delta(S)$. Since $x_0 \in X \setminus N_\delta(S)$, $x_0 \notin S_0$ and $S_0$ is maximal non-strictly $\delta$-separated set in $X \setminus N_\delta(S)$, then $d(x_0, x) < \delta$ for some $x \in S_0$ and it gives a contradiction.
\end{proof}

\label{S:DCtoAll}
\begin{theorem}[\textbf{ZFA}]
\label{T:DCImpliesEverything}
Let $\mathbf{\Psi}$ be one of $ \mathbf{PS, MS, PN, MN, PSE, MSE, PNE}$ or $\mathbf{MNE}$, then $\mathbf{DC} \implies \mathbf{\Psi}$.
\end{theorem}
\begin{proof}
Let $(\mathbf{\Psi_1}, \mathbf{\Psi_2})$ be one of (\textbf{PS}, \textbf{PSE}), (\textbf{MS}, \textbf{MSE}), (\textbf{PN}, \textbf{PNE}) or (\textbf{MN}, \textbf{MNE}). By Remark \ref{rem1} we have \textbf{DC} $\implies$ $\mathbf{\Psi_1}$. Thus by virtue of Proposition \ref{CCplusPsi} we get $\textbf{DC} \implies \textbf{DC} + \mathbf{\Psi_1} \implies \textbf{CC} + \mathbf{\Psi_1} \implies \mathbf{\Psi_2}$.
\end{proof}
\par  
  We now turn to a brief discussion of Fraenkel-Mostowski models and a transfer theorem of Pincus.  Fraenkel-Mostowski permutation models of the theory \textbf{ZFA} provide a relatively easy way of proving independence results in the theory \textbf{ZFA}.\footnote{The reader is referred to \cite[Chapter 4]{Jech} for a discussion of permutation models.}  
  These independence results usually have the form ``There is a permutation model in which $\mathbf{\Phi} = \mathbf{\Gamma} \land \lnot \mathbf{\Psi}$ holds'' where $\mathbf{\Gamma}$ and $\mathbf{\Psi}$ are consequences of \textbf{AC}.  A sentence $\mathbf{\Phi}$ is \emph{transferable} if there is a meta theorem  ``If there is a Franekel-Mostowski model of \textbf{ZFA} in which $\mathbf{\Phi}$ is true then there is a model of \textbf{ZF} in which $\mathbf{\Phi}$ is true''.  We will be using a transfer theorem of Pincus (\cite[Theorem 4]{Pincus}).  The part of the theorem that we need is
  \begin{theorem} \label{T:Pincus}
  Let $\mathbf{\Phi}$ be a conjunction of $\mathbf{CC}$ and any (finite) number of injectively boundable statements. If $\mathbf{\Phi}$ has a Fraenkel-Mostowski model then $\mathbf{\Phi}$ has a \textbf{ZF} model.  
  \end{theorem}
  For a definition of \emph{injectively boundable} see \cite{Pincus} or \cite[Note 103]{HR}.

\section{Strict $\delta$-separation in separable pseudometric spaces} \label{S:sepPsmSpaces}

We first show that $\textbf{PS}$ is provable in \textbf{ZFA}.
\begin{theorem}[\textbf{ZFA}]
\label{T:ZFtoPS}

For all $\delta \in \mathbb{R}^+$, $\PS$.
\end{theorem}
\begin{proof}
By Theorem \ref{T:DeltaDoesntMatter}, we only have to prove the Theorem for $\delta = 1$. As in the proof of Theorem \ref{T:DCImpliesEverything} this proof depends on ideas from \cite{Dyb-Gorka1}.  Let $(X,p)$ be a separable pseudometric space with countable dense subset $D = \{d_i : i \in \mathbb{Z}^{+} \}$ and assume that $(X,p)$ has no finite maximal strictly $1$-separated set.  

\begin{lemma} \label{L:ForDefOfkini}
If $T$ is a finite subset of $X$ and $N \in \mathbb{Z}^{+}$ is such that $\Omega := (d_N, 1/2) \setminus \bigcup_{t\in T}\overline{B}(t,1) \neq \emptyset$, then there exists $i \in \mathbb{Z}^{+}$ such that $d_i \in \Omega$.
\end{lemma}
\begin{proof}
The set $\Omega$ is open and by our assumption non-empty. Since $D$ is dense in $(X,p)$, there is an $i \in \mathbb{Z}^{+}$ such that $d_i \in \Omega$.
This completes the proof of the Lemma.
\end{proof}
We shall define by recursion a sequence of pairs of positive integers $(k_i, n_i)_{i=1}^\infty$ so that the set $\{ d_{k_i} : i \in \mathbb{Z}^+ \}$ is a maximal strictly $1$-separated set in $(X,p)$. We give the construction of the set $(k_i, n_i)_{i=1}^\infty$ and prove by induction that $(k_i, n_i)_{i=1}^\infty$ has the following Properties for $i \in \mathbb{Z}^+$.
\begin{enumerate}
\item \label{LI:dkInB} $d_{k_i} \in B(d_{n_i}, \frac{1}{2})$,
\item \label{LI:OrdPre} if $j \in \mathbb{Z}^+$ is such that $j < i$, then $n_j < n_i$,
\item \label{LI:pkdgt1} if $j \in \mathbb{Z}^+$ is such that $j < i$, then $p(d_{k_i}, d_{k_j}) > 1$.
\end{enumerate}
{\bf Construction}  of the sequence $(k_i, n_i)_{i=1}^\infty$:

 $\bullet (k_1,n_1) = (1,1)$.
It is clear that Properties (\ref{LI:dkInB}) through (\ref{LI:pkdgt1}) are true when $i = 1$.

$\bullet$ Assume that $(k_i, n_i)$ has been defined for all $i < m$ and that for all $i < m$, Properties (\ref{LI:dkInB}) through (\ref{LI:pkdgt1}) are true.  By Property (\ref{LI:pkdgt1}) the set $\{ d_{k_i} : i < m \}$ is strictly $1$-separated, but, by our assumption, is not maximal.  Therefore there exists $ y  \in X$ such that for all $i < m$, $p(y, d_{k_i}) > 1$ and thus there must be an $N \in \mathbb{Z}^+$ such that $y \in B(d_N, 1/2)$.  Let 
\begin{equation} \label{D:n_m}
 n_m = \min \left\{ l: B(d_l, 1/2) \setminus \bigcup_{i=1}^{m-1}\overline{B}(d_{k_i},1) \neq \emptyset \right\}. 
\end{equation}
  By Lemma \ref{L:ForDefOfkini} there is $k \in \mathbb{Z}^+$ such that $d_k \in B(d_{n_m}, 1/2) \setminus \bigcup_{i=1}^{m-1}\overline{B}(d_{k_i},1)$.  We let 
\begin{equation}  \label{D:k_m}
  k_m= \min \left\{ k: d_k \in B(d_{n_m}, 1/2) \setminus \bigcup_{i=1}^{m-1}\overline{B}(d_{k_i},1) \right\}.
\end{equation}  

Assuming that $(k_i, n_i)$ satisfies items (\ref{LI:dkInB}) through (\ref{LI:pkdgt1}) (preceding the Construction) for all $i < m$ we show that $(k_i, n_i)$ satisfies those properties for $i = m$.

Property (\ref{LI:dkInB}) follows from Equation (\ref{D:k_m}) in the Construction.

For Property (\ref{LI:OrdPre}) we first note that for all $j < m$, $n_j \ne n_m$.  Indeed, since $p(d_{n_j}, d_{k_j}) <1/2$ and $p(d_{n_m}, d_{k_m}) < 1/2$, the equality $n_j = n_m$ would lead to the contradiction 
\[
p(d_{k_j}, d_{k_m}) \le p(d_{k_j}, d_{n_j}) + p(d_{n_j}, d_{k_m}) = p(d_{k_j}, d_{n_j}) + p(d_{n_m}, d_{k_m}) < 1.
\]
 To complete the proof assume that Property (\ref{LI:OrdPre}) is false and let $i'$ be the least positive integer less than $m$ for which $n_m \le n_{i'}$.  Then, since for all $j <m$ $n_j \ne n_m$, we have $n_m < n_{i'}$.  We also know that $d_{k_m}\in  B(d_{n_m}, 1/2) \setminus \bigcup_{i=1}^{i'-1}\overline{B}(d_{k_i},1) $. Therefore, $ n_{i'} \le n_m$ and we get a contradiction.

Property (\ref{LI:pkdgt1}) follows from Equation (\ref{D:k_m}) in the Construction.

To complete the proof of Theorem \ref{T:ZFtoPS} we have to argue that $S = \{ d_{k_i} : i \in \mathbb{Z}^+ \}$ is a maximal strictly $1$-separated set.  The fact that $S$ is strictly $1$-separated follows from Property (\ref{LI:pkdgt1}). To prove that 
 $S$ is maximal we assume the contrary.  Then there exists $y \in X\setminus S$ such that for all $i \in \mathbb{Z}^{+}, p(y,d_{k_i}) > 1$.  There is an $n \in \mathbb{Z}^+$ such that $y \in B(d_n, 1/2)$ and we assume that $n$ is the least such positive integer.  Since $\{ n_i \}_{i \in \mathbb{Z}^+}$ is a strictly increasing sequence of integers there is the least $j \in \mathbb{Z}^+$ such that $n \le n_j$.  It is not possible that $n = n_j$ because, as above this leads to the contradiction 
\[
p(d_{k_j}, y) \le p(d_{k_j}, d_{n_j}) + p(d_{n_j}, y) = p(d_{k_j}, d_{n_j}) + p(d_{n}, y) < 1
\]
so $n < n_j$.  On the other hand $y \in B(d_n, 1/2) \setminus \bigcup_{i=1}^{j-1}\overline{B}(d_{k_i}, 1)$ and thus by the definition of $n_j$ in (\ref{D:n_m}) we have $n_j \le n$ which is an obvious contradiction.
\end{proof}

\begin{theorem}[\textbf{ZFA}]
\label{T:ExtensionOfDeltaSep}
$\mathbf{PSE} \implies \mathbf{CC}$.
\end{theorem}

\begin{proof}
By Theorem \ref{T:DeltaDoesntMatter}, it suffices to prove $\mathbf{PSE}(1) \implies \mathbf{CC}$.
\par
Assume that $\mathcal{Y} = \{Y_n : n \in \omega \}$ is a countable pairwise disjoint family of non-empty sets (and that the mapping $n \mapsto Y_n$ is one to one).  Our plan is to apply $\mathbf{PSE}(1)$ to a certain separable pseudometric space $(X,d)$ and a strictly $1$-separated set in $(X,d)$ to get a maximal strictly $1$-separated set.  From this we will be able to obtain a choice function for $\mathcal{Y}$.
\par
To construct $X$ we assume without loss of generality that 
\begin{equation}
\label{E:DjntnssCondition}
\left(\mathbb{R}^{(\mathbb{Z}^+)} \times \omega\right) \cap \left(\bigcup \mathcal{Y}\right) = \emptyset.
\end{equation}
For $i \in \mathbb{Z}^+$ we define two sequences $a_i = (a_i(k))_{k\in \mathbb{Z}^+}$ and $a_i' = ( a_i'(k))_{k \in \mathbb{Z}^+}$ (both in $\mathbb{R}^{(\mathbb{Z}^+)}$) by 
\begin{align*}
a_i(k) = \begin{cases}\frac{i+1}{i} & \mbox{ if } i = k\\
                       0 & \mbox{ otherwise}
         \end{cases} , \,\,\,\,
a_i'(k) = \begin{cases} \frac{1}{i} & \mbox{ if }i = k \\
                        0 & \mbox{ otherwise} .
          \end{cases}              
\end{align*}

Let $\mathbf{z} \in \mathbb{R}^{(\mathbb{Z}^+)}$ be the zero sequence defined by $\mathbf{z}(k) = 0$ for all $k \in \mathbb{Z}^+$ and define $A$ and $A_E$ respectively by
\begin{align*}
A = \{ a_i : i \in \mathbb{Z}^+ \} \cup \{ a_i' : i \in \mathbb{Z}^+ \}, \,\,\,\,
A_E = \{ a_i : i \in \mathbb{Z}^+ \} \cup \{ a_i' : i \in \mathbb{Z}^+ \} \cup \{ \mathbf{z} \}.
\end{align*}
Let $d_E$ be the usual Euclidean metric on $A_E$ defined by 
\[
d_E(b,c) = \sqrt{\sum_{k=1}^\infty (c(k)-b(k))^2}.
\]
For each $n \in \omega$, let $X_n = (A \times \{ n \}) \cup Y_n$ and let $X = \bigcup_{n \in \omega} X_n$. Finally, define a metric $d$ on $X$ by 
\[
d(t,s) = \begin{cases} 3 & \mbox{ if } \exists m,n \in \omega \mbox{ with } m \ne n \land  t \in X_n \land s \in X_m \\
                       0 & \mbox{ if } \exists n \in \omega \mbox{ such that }t,s \in Y_n \\
                       d_E(x,y) & \mbox{ if } \exists n \in \omega \hskip.05in \exists x, y \in A \mbox{ such that } t = (x,n) \land s = (y,n) \\
                       d_E(x,\mathbf{z}) & \mbox{ if } \exists n \in \omega \hskip.05in \exists x \in A \mbox{ such that } (t= (x,n) \land s \in Y_n ) \mbox{ or }\\
                       & \hskip.2in(s = (x,n) \land t \in Y_n)
         \end{cases}.
\]
Our assumption that the elements of $\mathcal{Y}$ are pairwise disjoint and the assumption (\ref{E:DjntnssCondition}) ensure that $d$ is well defined.
Furthermore, one can easily convince oneself that  $(X,d)$ is a pseudometric space. 
\par
We now list some properties of $(X,d)$ which will be useful in proving Lemmas \ref{L:Ddense} through \ref{L:SextGivesCF}. The proofs are straightforward.
\begin{lemma}\label{L:propsOfXd}
$\,$
\begin{enumerate}
\item \label{LI:da'ins}  For all $n \in \omega,  i \in \mathbb{Z}^+,  s \in Y_n, d((a'_i,n),s) = 1/i$. 
\item  \label{LI:daiaj>1} If $i \ne j$, $d((a_i,n),(a_j,n)) = \sqrt{\left( \frac{i+1}{i}\right)^2 + \left(\frac{j+1}{j}\right)^2} > 1$. 
\item \label{LI:daia'i=1} For all $ n \in \omega,  i \in \mathbb{Z}^+, d((a_i,n),(a'_i,n)) = 1$. 
\item \label{LI:dst<1tot=a'i} For all $n \in \omega, s \in Y_n, t \in X \setminus Y_n$, if $d(s,t) \leq 1$ then $t = (a'_i,n)$ for some $i \in \mathbb{Z}^+$. 
\end{enumerate}
\end{lemma}
\begin{lemma} \label{L:Ddense}
The set $D = \bigcup_{n \in \omega} (A \times \{ n \}) = A \times \omega$ is a countable dense subset of $(X,d)$.
\end{lemma}
\begin{proof}
Since $A$ is countable, $D = A \times \omega$ is countable. Furthermore, since $X = D \cup \left( \bigcup_{n \in \omega} Y_n \right)$, to show that $D$ is dense in $(X,d)$ it suffices to prove that every $t \in \bigcup_{n \in \omega} Y_n$ is a limit point of $D$.  By definition $\{ (a'_i,n)  : i \in \mathbb{Z}^+ \land n \in \omega \} \subseteq D$.  Therefore part (\ref{LI:da'ins}) of Lemma \ref{L:propsOfXd} implies $t$ is a limit point of $D$.
\end{proof}
Let $S = \bigcup_{n \in \omega} \left( \{a_i : i \in \mathbb{Z}^+ \} \times \{ n \} \right) = \{a_i : i \in \mathbb{Z}^+ \} \times \omega$. 
\begin{lemma} \label{L:SoneSeparated}
The set $S$ is strictly $1$-separated.
\end{lemma}
\begin{proof}
Assume $(a_i,n)$ and $(a_j,m)$ are distinct elements in $S$.  If $n \ne m$ then by the first clause of the definition of $d$, $d((a_i,n),(a_j,m)) = 3 > 1$.  On the other hand, if $n = m$ then $i \ne j$, so by part (\ref{LI:daiaj>1}) of Lemma \ref{L:propsOfXd}, $d((a_i,n),(a_j,m)) > 1$.
\end{proof}
\begin{lemma} \label{L:SextGivesCF}
If $S'$ is a maximal strictly $1$-separated set containing $S$ then for all $n \in \omega$, $|S' \cap Y_n| = 1$.
\end{lemma}
\begin{proof}
Assume $S'$ satisfies the hypotheses of the Lemma and assume $n \in \omega$.  Clearly, $|S' \cap Y_n| \le 1$.  We prove $|S' \cap Y_n| \ge 1$ by contradiction. Suppose $S' \cap Y_n = \emptyset$.  Since $S'$ is maximal, for all $y \in Y_n$ there is a $t \in S' \setminus Y_n$ such that $d(t,y) \le 1$. By part (\ref{LI:dst<1tot=a'i}) of Lemma \ref{L:propsOfXd},  $t = (a'_i,n)$ for some $i \in \mathbb{Z}^+$.  But by part (\ref{LI:daia'i=1}) in Lemma \ref{L:propsOfXd}, $d((a'_i,n), (a_i,n)) = 1$.  This is a contradiction since $(a_i,n) \in S \subseteq S'$.
\end{proof}
By Lemmas \ref{L:Ddense} and \ref{L:SoneSeparated} and $\mathbf{PSE}(1)$, there is a maximal strictly $1$-separated set $S'$ containing $S$.  Therefore, by Lemma \ref{L:SextGivesCF}, $\mathcal{Y}$ has a choice function.
\end{proof}

\section{Non-strict $\delta$-separation in separable pseudometric spaces} \label{S:PN}
\begin{theorem}[\textbf{ZFA}]
\label{T:NonStDelSepToCC}
$\mathbf{PN} \implies  \mathbf{CC}$.
\end{theorem}
\begin{proof}
As in the proof of Theorem \ref{T:ExtensionOfDeltaSep}, it suffices to prove $\mathbf{PN}(1) \implies \mathbf{CC}$.
Assume $\textbf{PN}(1)$ and let $X$ be a countable family of non-empty, pairwise disjoint sets where the function $n \mapsto X_n$ is one to one from $\omega$ onto $X$.  Also assume without loss of generality that $(\bigcup X) \cap \mathbb{R} = \emptyset$.  We will construct a choice function for $X$.
 \begin{lemma}
\label{L:delta=1}
There is a separable pseudometric space $(Z,p)$ such that
\begin{enumerate}[topsep=0pt]
\item \label{LC:subset} $\bigcup_{n \in \omega} (X_n \times \{ 0, 1 \}) \subseteq  Z$
\item \label{LC:cardnalities} Every maximal non-strictly $1$-separated set $S$ in $(Z,p)$ satisfies for all $n \in \omega$,
  \begin{enumerate}
  \item
   \label{LC:gt0}  
   $0 < \large| S \cap \large(X_n \times \{ 0, 1\}\large) \large| $
  \item 
   \label{LC:le1}  
  $\left|  S \cap (X_n \times \{ 0 \}) \right| \le 1$ and $\left|  S \cap (X_n \times \{ 1 \}) \right| \le 1$
  \end{enumerate}
\end{enumerate}
\end{lemma}
\par
Before giving the proof of the Lemma we argue that Lemma \ref{L:delta=1} is sufficient to prove Theorem \ref{T:NonStDelSepToCC} in the case $\delta = 1$.  
\par 
To construct a choice function $f$ for $X$ let $(Z,p)$ be the pseudometric given by Lemma \ref{L:delta=1}.  By $\textbf{PN}(1)$ there is a maximal non-strictly $1$-separated set $S$ in $(Z,p)$ and by Lemma \ref{L:delta=1}, for each $n \in \omega$, each of $S \cap (X_n \times \{ 0 \})$ and $S \cap (X_n \cap \{ 1 \})$ contains at most one element and  one of these two sets is non-empty.  We can therefore define $f(X_n)$ to be the element of $x \in X_n$ for which $(x,0) \in S$ if such an element exists.  Otherwise $f(X_n)$ is the element $x \in X_n$ for which $(x,1) \in S$.
\begin{proof}[Proof of Lemma \ref{L:delta=1}]  
We begin by defining a certain subset $Y_0$ of the Euclidean plane $(\mathbb{R}^2, d)$ as follows.  Let 
\begin{align*}
A_0^0 &= \left\{ (x,y) \in \mathbb{R}^2: x^2 + (y-1)^2 = 1 \mbox{ and } 0 < x < \sqrt{3}/2 \mbox{ and $x$ is rational and } 0 < y < 1/2 \right\},\\ A_0^1 &= \left\{ (x,y) \in \mathbb{R}^2 : x^2 + y^2 = 1 \mbox{ and } 0 < x < \sqrt{3}/2 \mbox{ and $x$ is rational and } 1/2 < y < 1 \right\}, \\ 
Y_0 &= A_0^0 \cup A_0^1.
\end{align*}
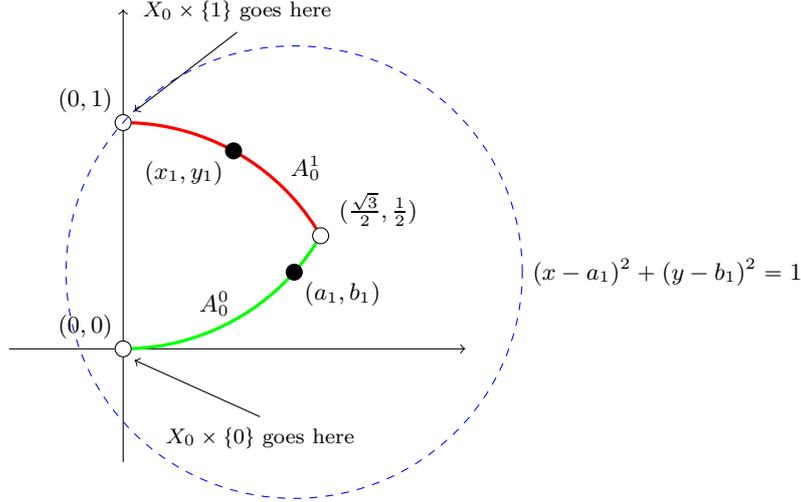
\begin{center}
 \begin{tikzpicture}[scale = 3]
 { 
 \draw[->] (-.5,0) -- (1.5,0) coordinate (x axis);
 \draw[->] (0,-.5) -- (0,1.5) coordinate (y axis);
 \draw[green, very thick, domain=-0:sqrt(3)/2] plot (\x, {-sqrt(1-\x^2)+1});
  \draw[red, very thick, domain=-0:sqrt(3)/2] plot (\x, {sqrt(1-(\x)^2)}); 
  \node[above right] at (.866,.5) {\footnotesize{ $(\frac{\sqrt{3}}{2},\frac{1}{2})$ }}; 
  \filldraw[white]   (.866,.5) circle (1.pt);
  \draw[black]   (.866,.5) circle (1.pt);
   \node[above left] at (0,0){\footnotesize{$(0,0)$}};
  \filldraw[white]   (0,0) circle (1.pt);
  \draw[black]   (0,0) circle (1.pt);   
  \node[above left] at (0,1){\footnotesize{$(0,1)$}};
  \filldraw[white]   (0,1) circle (1.pt);
  \draw[black]   (0,1) circle (1.pt);      
  }
  \node at (.8,.8) {\footnotesize{$A_0^1$}};
  \node at (.4,.2) {\footnotesize{$A_0^0$}};  
 \filldraw[black]   (.75,.339) circle (1.pt);     
 \node[below right] at (.75,.339){\footnotesize{$(a_1,b_1)$}};
 \filldraw[black]   (.4841,.875) circle (1.pt);     
 \node[below left] at (.4841,.875){\footnotesize{$(x_1,y_1)$}}; 
 \draw[blue, dashed] (.75,.339) circle (1);  
 \node[right] at (1.75,.339){\footnotesize{$(x-a_1)^2 + (y - b_1)^2 =  1$}}; 
  \draw[<-] (0.05,1.05) -- (.5,1.4) node[above]{\scriptsize  $X_0 \times \{ 1 \}$ goes here};
    \draw[<-] (0.05,-0.05) -- (.6,-.3) node[below]{\scriptsize  $X_0 \times \{0\}$ goes here};
\end{tikzpicture}
\end{center}
Note that both $A_0^0$ and $A_0^1$ are countable.  For every positive integer $n$ we let $A_n^0 = \{ (a, b + 2n) : (a,b) \in A_0^0 \}$, $A_n^1  = \{ (x, y + 2n) : (x,y) \in A_0^1 \}$ and  $Y_n = A_n^0 \cup A_n^1$. Since for every $(x,y) \in Y_0$, we have $0 < y < 1$ we conclude that for every pair $(x,y) \in Y_n$ we have $2n < y < 2n + 1$.
\begin{sublemma}
\label{L:A_0A_1close}
Let $n \in \omega$, then for all $(a_1,b_1)$ and $(a_2, b_2)$ in $A_n^0$ and for all $(x_1,y_1)$ and $(x_2,y_2)$ in $A_n^1$ we have
\begin{align}
\label{E:PropsOfA_n1}
d((a_1,b_1),(x_1,  y_1)) < 1, & \\ 
\label{E:PropsOfA_n2}
d((a_1,b_1), (a_2,b_2)) < 1 &\mbox{ and } d((x_1,y_1), (x_2,y_2)) < 1, \\ 
\label{E:PropsOfA_n3}
d((a_1,b_1),(0,2n)) < 1 &\mbox{ and } d((x_1,y_1),(0,2n+1)) < 1, \\
d(a_1,b_1), (0,2n+1)) = 1 &\mbox{ and } d((x_1,y_1),(0,2n)) = 1.
\label{E:PropsOfA_n4}
\end{align} 
\end{sublemma}
\begin{proof}
It is enough to prove the Sublemma for $n=0$. Parts (\ref{E:PropsOfA_n2}), (\ref{E:PropsOfA_n3}) and (\ref{E:PropsOfA_n4}) are clear from the picture. For part (\ref{E:PropsOfA_n1}) an algebraic proof can be given by assuming that $(a_1,b_1) \in A_0^0$ and $(x_1,y_1) \in A_0^1$ and calculating the square of the distance from $(a_1,b_1)$ to $(x_1,y_1)$.
\end{proof}

For $n \in \omega$, let $Z_n = Y_n \cup \left( X_n \times \{0, 1\}\right)$ and let $Z = \bigcup_{n\in \omega} Z_n$.  We shall define a pseudometric $p$ on $Z$ using the usual metric on $\mathbb{R}^2$ and for each $n \in \omega$, placing one copy of $X_n$ (namely $X_n \times \{ 0 \}$) at the point $(0,2n)$ and one copy of $X_n$ (namely $X_n \times \{ 1 \}$)  at $(0,2n+1)$.  Here are the details. Define a function $c : Z \to \mathbb{R}^2$ by
\[
c(t,s) = \begin{cases}  (t,s) & \mbox{ if } (t, s) \in \bigcup_{n \in \omega} Y_n \\
(0, 2n) & \mbox{ if } t \in X_n \mbox{ and } s = 0 \\
           (0,2n+1) & \mbox{ if } t \in X_n \mbox{ and } s = 1  \end{cases}.
\]
Function $c$ assigns each point of $Z$ to a point in the plane and the distance between two elements of $Z$ will be the Euclidean distance between their assignments.  That is, define a pseudometric $p$ on $Z$ by 
\begin{equation}
\label{D:defOfp}
p((t_1,s_1),(t_2,s_2)) = d(c(t_1,s_1),c(t_2,s_2)).
\end{equation}

\begin{sublemma}
\label{SL:propsOfp}
Assume $z_1, z_2 \in A_n^0$, $w_1, w_2 \in A_n^1$, $u_1, u_2 \in X_n \times \{ 0 \}$ and $v_1, v_2 \in X_n \times \{ 1 \}$, then
\begin{enumerate} 
\item \label{SLI:no1} $p(z_1, w_1) < 1$,
\item $p(z_1, z_2) < 1$ and $p(w_1,w_2) < 1$,
\item $p(z_1,u_1) < 1$ and $p(w_1,v_1) < 1$,
\item \label{SLI:No4} $p(z_1, v_1) = 1$ and $p(w_1,u_1) = 1$,
\item \label{SLI:No5} $p(u_1,v_1) = 1$,
\item \label{SLI:deq0} $p(u_1,u_2) = 0$ and $p(v_1,v_2) = 0$.
\end{enumerate}
\end{sublemma}
\begin{proof}
Using the definitions of $c$ and $p$, items (1) through (4) follow from equations (\ref{E:PropsOfA_n1}) through (\ref{E:PropsOfA_n4}) respectively and items (5) and (6) are clear from the definitions of $c$ and $p$.  
\end{proof}

It is clear that $(Z,p)$ is a pseudometric space. It is also the case that $\bigcup_{n \in \omega} Y_n$ is a countable dense subset of $Z$ since every element of $X_n \times \{ 0 \}$ is a limit point of $A_n^0$ and every element of $X_n \times \{ 1 \}$ is a limit point of $A_n^1$.  Therefore $(Z,p)$ is a separable pseudometric space satisfying condition (\ref{LC:subset}) of Lemma \ref{L:delta=1}.  
\par
To argue that condition (\ref{LC:cardnalities}) of Lemma \ref{L:delta=1} holds assume that $S$ is a maximal non-strictly $1$-separated set in $(Z,p)$ and that $n \in \omega$.  Part (\ref{LC:le1}) follows from Sublemma \ref{SL:propsOfp} part (\ref{SLI:deq0}).  To argue for part (\ref{LC:gt0}), we begin by showing that for $i \in \{0,1\}$
\begin{equation}
\label{E:Ifpuz<1}
\forall u \in X_n \times \{ i \}, \forall z \in Z \setminus (X_n \times \{ i \})( p(u,z) < 1 \Rightarrow z \in A^i_n ).
\end{equation}

We will argue for $i=0$ and the case $i=1$ we leave to the reader. Assume $u \in X_n \times \{ 0 \}$, $z \in Z \setminus (X_n \times \{ 0 \})$ and $p(u,z) < 1$.  Under these assumptions $z \notin Z_m$ where $m \in \omega \setminus \{ n \}$ because for every element $z$ of $Z_m$, if $r_1$ is the second component of $c(z)$ then $2m \le r_1 \le 2m + 1$.  Whereas, if $r_2$ is the second component $c(u)$ then $r_2 = 2n$.  This implies that $p(u,z) = d(c(u),c(z)) \ge 1$.  Since $z \notin \bigcup_{m \in \omega \setminus \{ n \}} Z_m$ we have $z \in Z_n \setminus (X_n \times \{ 0 \}) =  A^0_n \cup A^1_n \cup  (X_n \times \{ 1 \})$.  
By Sublemma \ref{SL:propsOfp}, item (\ref{SLI:No4}), $z \notin A^1_n$.  By Sublemma \ref{SL:propsOfp}, item (\ref{SLI:No5}), $z \notin X_n \times \{ 1 \}$.  
Therefore $z \in A^0_1$ completing the proof of (\ref{E:Ifpuz<1}).
\par 
We complete the proof of Lemma \ref{L:delta=1}, part (\ref{LC:cardnalities}) by contradiction.  Assume that both $S \cap (X_n  \times \{ 0 \})$ and $S \cap (X_n \times \{ 1 \})$ are empty.  By the first of these assumptions, since $S$ is maximal, for every $z \in X_n \times \{ 0 \}$ there is an element $s \in S$ such that $p(z,s) < 1$. Hence, by equation (\ref{E:Ifpuz<1}) this means that there is an element $e_0 \in S \cap A^0_n$.  Similarly, there is an element $e_1 \in S \cap A^1_n$.  But this is a contradiction, since by Sublemma \ref{SL:propsOfp}, part (\ref{SLI:no1}), $p(e_0,e_1) < 1$.
\end{proof}
\end{proof}
\begin{corollary} \label{C:CCiffPSE}
\begin{enumerate}
\item[]
\item[(i)] $\textbf{CC} \iff \textbf{PSE} \implies \textbf{MSE}$;
\item[(ii)] $ \textbf{PN} \iff \textbf{PNE} \implies \textbf{CC}$.
\end{enumerate}
    
\end{corollary}

\begin{proof}
With the use of Proposition \ref{CCplusPsi} statement (i) follows from  Theorems \ref{T:ZFtoPS} and \ref{T:ExtensionOfDeltaSep} and statement (ii) follows from  Theorem \ref{T:NonStDelSepToCC}. 
\end{proof}

\section{Strict and non-strict $\delta$-separation in separable metric spaces}\label{S:InMetricSpaces}

We first make several remarks about the consequences of Proposition \ref{P:XleR} for  Fraenkel-Mostowski models of \textbf{ZFA}.
\begin{remark} \label{R:FMmodels} \leavevmode
\begin{enumerate}
\item \label{I:M_InFM} In any Fraenkel-Mostowski model of \textbf{ZFA} the real numbers are well orderable and therefore, by Proposition \ref{P:XleR}, every separable metric space is well orderable.  It follows from the discussion in \cite[Section 4.2]{Jech} that all of the statements $\mathbf{MSE}$, $\mathbf{MNE}$, $\mathbf{MS}$ and $\mathbf{MN}$ are true in all Fraenkel-Mostowski models.
\item \label{I:M_Trans} The statements $\mathbf{MSE}$, $\mathbf{MNE}$, $\mathbf{MS}$ and $\mathbf{MN}$ are injectively boundable.  See \cite[Note 103]{HR} or \cite{Pincus} for a definition of \emph{injectively boundable}.
\item Therefore using Theorem \ref{T:Pincus}, if $\mathbf{\Phi}$ is any consequence of $\mathbf{AC}$ which is false in some Fraenkel-Mostowski model and whose negation is injectively boundable then none of the statements $\mathbf{MSE}$, $\mathbf{MS}$, $\mathbf{MNE}$ or $\mathbf{MN}$ implies $\mathbf{\Phi}$ in $\mathbf{ZF}$.  For most of the sentences $\mathbf{\Phi}$ which are consequences of \textbf{AC} and mentioned in \cite{HR}, the negation of $\mathbf{\Phi}$ is injectively boundable.  In particular the negations of $\mathbf{CC}$ and $\mathbf{DC}$ are injectively boundable.  So none of these four statements implies $\mathbf{CC}$ nor does any one of them imply $\mathbf{DC}$.  See \cite[Theorem 4]{PincusADC} or Note 103 of \cite{HR} for a more complete discussion of the transfer theorem of Pincus used here.
\end{enumerate}
\end{remark}

\begin{theorem}[\textbf{ZFA}]\label{MSEonR}
$\mathbf{CC}(\mathbb{R}) \iff \mathbf{MSE}$.
\end{theorem}

\begin{proof}
    $\mathbf{MSE} \implies \textbf{CC}(\mathbb{R})$: Assume that $\mathcal{Y} = \left\{Y_n \colon n \in \omega\right\}$ is a countable family of pairwise disjoint non-empty subsets of $\mathbb{R}$ and that the mapping $n \mapsto Y_n$ is one-to-one. Since $\mathbb{R}$ is bijective with $(3n, 3n+1)$ for all $n \in \omega$, we can assume that $Y_n \subset (3n, 3n+1)$ for all $n \in \omega$. Let $Y = \bigcup \mathcal{Y}$. Let $D_0$ be the set of dyadic numbers i.e. $D_0 = \left\{m/ 2^n \colon\, m \in \mathbb{Z}, n \in \mathbb{Z}^+\right\}$ and $D_1 = D_0  + 1/3$. Sets $D_0$ and $D_1$ are countable and dense in $\mathbb{R}$ and it is easy to check that $D_0 \cap D_1 = \emptyset$. We can assume\footnote{If $Y_n \cap D_i \neq \emptyset$ for some $n \in \omega$ and $i=0, 1$, then we can take $j_0 = \min\left\{j \in \mathbb{Z}^+\colon\; d^i_j \in Y_n \cap D_i\right\}$ where $D_i = \left\{d_j^i\right\}_{j=1}^\infty$. Thus we can put $y_n = d_{j_0}^i \in Y_n$. Hence from those sets $Y_n$ such that $Y_n \cap \brac{D_0 \cup D_1} \neq \emptyset$ we can select elements directly, without any choice.} that $Y \cap \brac{D_0 \cup D_1} = \emptyset$. Finally, we put $X = Y \cup D_0 \cup D_1$ and define metric $d$ on $X$ in the following manner:

 \[
d(x, y) = \begin{cases}
0 & \text{if } x=y\\
1 + |x - y| & \text{if } x\neq y \; \text{and} \; x,y \in D_1 \\
1 + |y - x + 1/3| & \text{if } x \in D_1 \; \text{and} \; y \not\in D_1 \\
1 + |x - y + 1/3| & \text{if } x \not\in D_1 \; \text{and} \; y \in D_1 \\
|x-y| & \text{otherwise}
\end{cases}.
\]
It is easy to check that $d$ is indeed metric and $(X, d)$ is separable metric space with $D_0 \cup D_1$ as the countable and dense set in $X$. Since $D_1$ is strictly $1$-separated set in $X$ and $X \subset \mathbb{R}$ is separable, we can use the assumed statement so there exists a maximal strictly $1$-separated set $S'$ in $X$ such that $D_1 \subset S'$. We shall show that $|Y_n \cap S'| = 1$ for all $n \in \omega$.

Let us notice that for every $x \in D_0$ there exists $y \in D_1$ such that $d(x, y) = 1$. Indeed, it is enough to take $y = x + 1/3 \in D_1$ and then $d(x, y) = 1 + |x - y + 1/3| = 1$. Since $S'$ is strictly $1$-separated set and $D_1 \subset S'$, it means that $S' \cap D_0 = \emptyset$. 

Let us fix $n \in \omega$. If $x, y \in Y_n \cap S' \subset Y_n \subset (3n, 3n + 1)$, then $d(x, y) = |x - y| < 1$. Thus $x = y$ since $x, y \in S'$. Hence $|Y_n \cap S'| \le 1$.

Let us suppose that $Y_n \cap S' = \emptyset$ for some $n \in \omega$. Let us take $x \in Y_n$ and fix $y \in S'$. Since $S' \cap D_0 = \emptyset$, either $y \in D_1$ or $y \in Y$. If $y \in D_1$, then $y - 1/3 \in D_0$ and $Y_n \cap D_0 = \emptyset$ so $x \neq y - 1/3$. Thus $d(x, y) = 1 + |x - y + 1/3| > 1$. If $y \in Y$, then $y \not\in Y_n$ since $Y_n \cap S' = \emptyset$ so $d(x, y) > 1$. We have $d(x, y) > 1$ for all $y \in S'$ so it contradicts the maximality of $S'$. Hence $Y_n \cap S' \neq \emptyset$ for all $n \in \omega$.

$\textbf{CC}(\mathbb{R}) \implies \mathbf{MSE}$: By Proposition \ref{CCplusPsi} we have that $\textbf{CC}(\mathbb{R}) + \textbf{MS} \implies \mathbf{MSE}$ but \textbf{MS} is true in \textbf{ZFA} by  Theorem \ref{T:ZFtoPS}.
\end{proof}

\begin{corollary}\label{MSEonRcor}
\textbf{MSE} is not provable in \textbf{ZF}.
\end{corollary}

\begin{proof}
    It follows immediately from Theorem \ref{MSEonR}.
\end{proof}

\begin{proposition}[\textbf{ZFA}]\label{MNEimpliesCCR}
$\textbf{MNE} \implies \textbf{CC}(\mathbb{R})$. 
\end{proposition}

\begin{proof}
Assume that $\mathcal{Y} = \left\{Y_n \colon n \in \omega\right\}$ is a countable family of pairwise disjoint non-empty subsets of $\mathbb{R}$ and that the mapping $n \mapsto Y_n$ is one-to-one. Let $S^1$ be unit circle in $\mathbb{R}^2$, namely $S^1 = \left\{(x, y) \in \mathbb{R}^2\colon\; x^2 + y^2 = 1\right\}$. There exists a bijection $f\colon \mathbb{R} \to S^1$. Let $\widehat{Y}_n = f(Y_n) + (4n, 0)$ for all $n \in \omega$. Since $f$ is bijective and $\mathbb{R}$ is linearly ordered, it is enough to prove that there exists set $S'$ such that $S' \cap \widehat{Y}_n$ is non-empty and finite for all $n \in \omega$ to complete the proof.

Let $D^1 = \left\{(x, y) \in \mathbb{Q}^2\colon\; x^2 + y^2 < 1\right\}$ and $D_n^1 = D^1 + (4n, 0)$ for all $n \in \omega$. We consider metric space $(X, d)$ where $X = \bigcup_{n \in \omega} \brac{\widehat{Y}_n \cup D_n^1}$ and $d$ is the usual Euclidean metric on $\mathbb{R}^2$. Obviously, $X$ is separable with $\bigcup_{n \in \omega} D_n^1$ as a countable and dense set. Let $S = \left\{(4n, 0)\colon\; n \in \omega\right\}$. It is non-strictly $1$-separated set in $X$ so by \textbf{MNE} there exists a maximal non-strictly $1$-separated set $S'$ such that $S \subset S'$.

Let us notice that for all $n \in \omega, x \in D_n^1$ and $y \in \widehat{Y}_n$ we have $d\brac{(4n, 0), x} < 1$ and $d\brac{(4n, 0), y} = 1$. Thus $S' \cap \bigcup_{n \in \omega} D_n^1 \setminus \{(4n,0)\} = \emptyset$. Since $S'$ is maximal and $d(x, y) \ge 1$ for all $x \in D_n^1 \cup \widehat{Y}_n$ and $y \in D_m^1 \cup \widehat{Y}_m$ where $n \neq m$, it is easy to see that $S' \cap \widehat{Y}_n \neq \emptyset$ for all $n \in \omega$. Moreover for all $n \in \omega$ set $S' \cap \widehat{Y}_n$ is finite since every $1$-separated set contained in $S^1 + (4n, 0)$ is finite.
\end{proof}

\begin{corollary}[\textbf{ZFA}]
$\textbf{CC}(\mathbb{R}) + \textbf{MN} \iff \textbf{MNE}$.
\end{corollary}

\begin{proof}
It follows from the Proposition \ref{MNEimpliesCCR} and Proposition \ref{CCplusPsi}.
\end{proof}

\begin{theorem}[\textbf{ZFA}] \label{T:MN+CCtoPN}
$\textbf{MN} + \textbf{CC} \iff \textbf{PN}$.
\end{theorem}

\begin{proof}
For the implication ``$\Longleftarrow$'' we note that $\mathbf{PN} \implies \mathbf{MN}$ is clear and from Theorem \ref{T:NonStDelSepToCC} we have $\mathbf{PN} \implies \mathbf{MN}$.

For the other implication let $(X, p)$ be a separable pseudometric space and $D$ be a countable and dense set. We convert this space into metric space in the standard way. We define an equivalence relation $\sim$ on $X$ as follows: $x \sim y \Longleftrightarrow p(x,y) = 0$ for all $x, y \in X$. Then $(X/{\sim}, d)$ is metric space where $d\brac{[x], [y]} = p(x,y)$ for all $x, y \in X$ and $[x]$ denotes the equivalence class of $x$. Obviously $D/{\sim}$ is a countable and dense set in $X/{\sim}$ so $(X/{\sim}, d)$ is separable.

Since we assumed \textbf{MN}, there exists a maximal non-strictly $1$-separable set $S$ in $(X/{\sim}, d)$. Every non-strictly $1$-separated in every separable metric space is at most countable so $S$ is at most countable. Let $S = \left\{S_i\right\}_{i=1}^\infty$. From \textbf{CC} there exists set $\widehat{S} = \left\{s_i\right\}_{i=1}^\infty \subset X$ such that $s_i \in S_i$ for every $i \in \mathbb{Z}^+$. Then $[s_i] = S_i$ for every $i \in \mathbb{Z}^+$ and it is easy to see that $\widehat{S}$ is a maximal non-strictly $1$-separated set in $(X, p)$.
\end{proof}

\begin{corollary} \label{C:PNnotToDC}
In the theory \textbf{ZF}, the sentence \textbf{PN} (and therefore the sentence \textbf{PNE}) does not imply \textbf{DC}.
\end{corollary}
\begin{proof}
In \cite{HR1996} a Fraenkel-Mostowski model is constructed in which \textbf{CC} is true and \textbf{DC} is false. By Remark \ref{R:FMmodels} Part (\ref{I:M_InFM}), \textbf{MN} is also true in this model.  By Remark \ref{R:FMmodels} Part (\ref{I:M_Trans}), \textbf{MN} is injectively boundable.  It follows from Theorem \ref{T:Pincus} that there is a model of \textbf{ZF} in which both \textbf{CC} and \textbf{MN} are true and \textbf{DC} is false. 
Therefore, by Theorem \ref{T:MN+CCtoPN}, in this model \textbf{PN} is true and \textbf{DC} is false.
\end{proof}

We close this section with the following observation.
\begin{theorem}[\textbf{ZFA}]\label{nonstrictZF}
 Let $(X, d)$ be a separable, complete metric space and let $\delta > 0$. Then, for every non-strictly $\delta$-separated subset $S$ of $X$ there exists a maximal non-strictly $\delta$-separated subset $S'$ of $X$ such that $S \subset S'$.
\end{theorem}

\begin{proof}

First of all we shall prove the following lemma.
\begin{lemma}[\textbf{ZFA}]\label{choiceclosedcomplete}
Let $(X, d)$ be a separable, complete metric space and $\mathcal{F}$ be the family of all non-empty closed subsets of $X$. Then, there exists choice function for $\mathcal{F}$ i.e. set $\left\{a_F\right\}_{F \in \mathcal{F}} \subset X$ such that $a_F \in F$ for all $F \in \mathcal{F}$.
\end{lemma}

\begin{proof}
Let $\left\{x_i\right\}_{i=1}^\infty \subset X$ be a countable and dense subset of $X$. Let $F \in \mathcal{F}$. We shall give the construction of $a_F \in F$. We define sequence $\left\{n_k\right\}_{k=1}^\infty \subset \omega$ in the following manner: 
\begin{itemize}
\item $n_1 = \min\left\{i \in \mathbb{Z}^+\colon\; x_i \in \bigcup_{x \in F} B\brac{x, 1/2}\right\}$,
\item $n_{k+1} = \min\left\{i \in \mathbb{Z}^+\colon\; x_i \in B\brac{x_{n_k}, 1/2^k} \cap \bigcup_{x \in F} B(x, 1/2^{k+1})\right\}$ for all $k \ge 1$.
\end{itemize}

Obviously, set $\bigcup_{x \in F} B\brac{x, 1/2}$ is open and non-empty so $n_1$ is well-defined. Let us assume that $n_k$ are defined for $k \le m$. We will show that $n_{m+1}$ is well-defined. Since set $B\brac{x_{n_m}, 1/2^m} \cap \bigcup_{x \in F} B(x, 1/2^{m+1})$ is open, it is enough to show that it is non-empty. Let us notice that $x_{n_m} \in \bigcup_{x \in F} B\brac{x, 1/2^m}$ so $B\brac{x_{n_m}, 1/2^m} \cap F \neq \emptyset$. As a consequence $B\brac{x_{n_m}, 1/2^m} \cap \bigcup_{x \in F} B(x, 1/2^{m+1}) \neq \emptyset$ and it shows that $n_{m+1}$ is well-defined.
The sequence $\left\{n_k\right\}_{k=1}^\infty$ has the following properties:
\begin{enumerate}
\item[(1)] $x_{n_{k+1}} \in B\brac{x_{n_k}, 1/2^{k}}$ for all $k \ge 1$,
\item[(2)] $\text{dist}\brac{x_{n_k}, F} < 1/2^k$ for all $k \ge 1$.
\end{enumerate}

From $(1)$ we conclude that sequence $\left\{x_{n_k}\right\}_{k=1}^\infty$ is Cauchy sequence so it converges since $X$ is complete. Let $\lim_{k\to\infty} x_{n_k} = x \in X$. From $(2)$, since $\text{dist}\brac{x, F} \le \text{dist}\brac{x_{n_k}, F} + d(x, x_{n_k})$ for every $k \ge 1$, then $\text{dist}\brac{x, F} = 0$. But $F$ is closed so $x \in F$. We put $a_F := x$.
\end{proof}
Let $\left\{x_i\right\}_{i=1}^\infty \subset X$ be a countable and dense subset of $X$ and let $S$ be non-strictly $\delta$-separated subset of $X$. By virtue of Lemma \ref{choiceclosedcomplete} let $\mathcal{F}$ be the family of all non-empty closed subset of $X$ and $\left\{a_F\right\}_{F \in \mathcal{F}} \subset X$ be a set such that $a_F \in F$ for all $F \in \mathcal{F}$. Let us denote $N_\delta(S) = \bigcup_{x \in S} B\brac{x, \delta}$. If $N_\delta(S) = X$, then $S$ is maximal non-strictly $\delta$-separated set so we assume $N_\delta(S) \neq X$. Let $y \in X \setminus N_\delta(S)$. We define sequences $\left\{y_n\right\}_{n=1}^\infty \subset X$ and $\left\{k_n\right\}_{n=1}^\infty \subset \omega $ as follows\footnote{It may happen that construction of these sequences stops at some point i.e. we constructed $k_n$ for $n \le m$ and $\overline{B}\brac{x_i, \delta/3} \setminus \brac{\bigcup_{j=1}^{m} B\brac{y_j, \delta} \cup N_\delta(S)} = \emptyset$ for all $i > k_m$. The proof in that case is similar as in the case where construction of $k_n$ goes to infinity.}:

\begin{itemize}
\item $k_1 = 0, y_1 = y$,
\item $k_n = \min\left\{i > k_{n-1}\colon\; \overline{B}\brac{x_i, \delta/3} \setminus \brac{\bigcup_{j=1}^{n-1} B\brac{y_j, \delta} \cup N_\delta(S)} \neq \emptyset\right\}$ for all $n \ge 2$,
\item $y_n = a_F$ where $F = \overline{B}\brac{x_{k_n}, \delta/3} \setminus \brac{\bigcup_{j=1}^{n-1} B\brac{y_j, \delta} \cup N_\delta(S)}$ for all $n \ge 2$.
\end{itemize}

We shall show that $S' := S \cup \left\{y_n\right\}_{n=1}^\infty$ is a maximal non-strictly $\delta$-separated set in $X$ containing $S$. Let us fix $n \in \mathbb{Z}^+$. Since $y_n \notin \brac{\bigcup_{j=1}^{n-1} B\brac{y_j, \delta} \cup N_\delta(S)}$, then $d(y_n, y_j) \ge \delta$ for all $j < n$ and $d(y_n, s) \ge \delta$ for all $s \in S$. But $n \in \mathbb{Z}^+$ was arbitrary so $d(y_n, y_j) \ge \delta$ for all $n, j \in \mathbb{Z}^+$, $n \neq j$. It proves that $S'$ is non-strictly $\delta$-separated set.

Let us suppose that $S'$ is not maximal i.e. there exists $x \in X$ such that $d(x, y_n) \ge \delta$ for all $n \in \mathbb{Z}^+$ and $d(x, s) \ge \delta$ for all $s \in S$. There exists $i \in \mathbb{Z}^+$ such that $x \in \overline{B}\brac{x_i, \delta/3}$. Thus $x \in \overline{B}\brac{x_i, \delta/3} \setminus \brac{\bigcup_{j=1}^{\infty} B\brac{y_j, \delta} \cup N_\delta(S)}$. Let us suppose that $i \neq k_n$ for any $n \in \mathbb{Z}^+$. Sequence $\left\{k_n\right\}_{n=1}^\infty$ is increasing so there exists $n \in \mathbb{Z}^+$ such that $k_n < i < k_{n+1}$. By the definition of $k_{n+1}$ we have $\overline{B}\brac{x_i, \delta/3} \setminus \brac{\bigcup_{j=1}^{n} B\brac{y_j, \delta} \cup N_\delta(S)} = \emptyset$ so $\overline{B}\brac{x_i, \delta/3} \setminus \brac{\bigcup_{j=1}^{\infty} B\brac{y_j, \delta} \cup N_\delta(S)} = \emptyset$ which contradicts the fact that $x$ is an element of this set. Hence $i = k_n$ for some $n \in \mathbb{Z}^+$. Obviously $n \ge 2$ since $i \in \mathbb{Z}^+$ and $k_1 = 0 \notin \mathbb{Z}^+$. It means that $x \in \overline{B}\brac{x_{k_n}, \delta/3} \setminus \brac{\bigcup_{j=1}^{n-1} B\brac{y_j, \delta} \cup N_\delta(S)}$. Since $y_n \in \overline{B}\brac{x_{k_n}, \delta/3} \setminus \brac{\bigcup_{j=1}^{n-1} B\brac{y_j, \delta} \cup N_\delta(S)}$ either, we obtain $x, y_n \in \overline{B}\brac{x_{k_n}, \delta/3}$ so $d(x, y_n) \le 2\delta /3$. It gives a contradiction with the fact that $d(x, y_n) \ge \delta$. Hence $S'$ is maximal non-strictly $\delta$-separated set.
\end{proof}   

\bibliographystyle{amsalpha}

\begin{thebibliography}{99}

\bibitem{bentley} H. L. Bentley and H. Herrlich, \textit{Countable choice and pseudometric spaces}, Topology and its Applications,
\textbf{85}, 153-164, 1998.

\bibitem{Dyb-Gorka1} M. Dybowski and P. G\'{o}rka, \textit{The axiom of choice in metric measure spaces and maximal $\delta$-separated sets}, Archive for Mathematical Logic \textbf{62}, 735-749, 2023.

\bibitem{Gorka1} P. G\'{o}rka, \textit{Separability of a metric space is equivalent to the existence of a Borel measure}, The American Mathematical Monthly, \textbf{128}, 84--86, 2020.

\bibitem{herrlich} H. Herrlich, \textit{Axiom of choice}, Springer, 2006.

\bibitem{HR} P. Howard and J. E. Rubin, \textit{Consequences of the Axiom of Choice}, Mathematical Surveys and Monographs (\textbf{59}),  American Mathematical Society, Providence, RI, 1998.

\bibitem{HR1996} P. Howard and J. E. Rubin, \textit{The Boolean prime ideal theorem plus countable choice do not imply dependent choice}, Math\. Logic Quart., {\bf 42}, 410-420, 1996.

\bibitem{Jech} T. J. Jech, \textit{The Axiom of Choice}, Studies in Logic and the Foundations of Mathematics \textbf{75}, North-Holland Publising Co., Amsterdam, 1973.

\bibitem{Jensen} R. B. Jensen, \textit{Consistency results for \textbf{ZF}}, Notices  Am. Math. Soc. \textbf{14}, 137, 1967.

\bibitem{Pincus} D. Pincus, \textit{Zermelo-Fraenkel consistency results by Fraenkel-Mostowski methods}, J. Symbolic Logic \textbf{37}, 721--743, 1972.

\bibitem{PincusADC} D. Pincus, \textit{Adding dependent choice}, Ann. Math. Logic \textbf{11}, 105--145, 1977.

\bibitem{Tach1} E. Tachtsis,  \textit{On the existence of almost disjoint and MAD familits without \textsf{AC}},  Bulletin Polish Acad. Sci. Math. \textbf{67}, 101-124, 2019.


\end{thebibliography}

\section*{Conflict of interest}
 The authors declare that they have no conflict of interest.

\end{document}